\documentclass[11pt]{article}
\usepackage[T1]{fontenc}
\usepackage{comment}
\usepackage[utf8]{inputenc}
\usepackage{amsmath,amssymb,amsthm,mathtools}
\usepackage[margin=1in]{geometry}
\usepackage{microtype,booktabs}
\usepackage{authblk}
\usepackage[numbers,sort&compress]{natbib}
\usepackage[colorlinks=true,linkcolor=blue,citecolor=blue,urlcolor=blue]{hyperref}
\allowdisplaybreaks
\DeclareMathOperator{\ex}{ex}
\newcommand{\cF}{\mathcal F}
\newcommand{\cW}{\mathcal W}
\newtheorem{theorem}{Theorem}[section]
\newtheorem{lemma}[theorem]{Lemma}
\newtheorem{proposition}[theorem]{Proposition}
\newtheorem{corollary}[theorem]{Corollary}
\newtheorem{question}[theorem]{Question}
\theoremstyle{definition}
\newtheorem{definition}[theorem]{Definition}
\theoremstyle{remark}
\newtheorem{remark}[theorem]{Remark}
\newtheorem{example}[theorem]{Example}

\title{Set-coloring Ramsey numbers with multiple intersections}
\author[1]{Guantao Chen}
\author[2]{Laihao Ding}
\author[2]{Qiong Fan\thanks{Corresponding author: \texttt{xfanqiong@ccnu.edu.cn}.}}
\author[2]{Xiaolan Hu}
\author[3]{Suyun Jiang}
\author[2]{Chunxiang Wang}
\affil[1]{Department of Mathematics and Statistics, Georgia State University, Atlanta, GA 30303, USA}
\affil[2]{School of Mathematics and Statistics, Key Laboratory of Nonlinear Analysis \& Applications (Ministry of Education), Central China Normal University, Wuhan, China}
\affil[3]{School of Artificial Intelligence, Jianghan University, Wuhan, Hubei, China}

\date{}
\hypersetup{pdftitle={Set-coloring Ramsey numbers with multiple intersections},pdfauthor={Guantao Chen, Laihao Ding, Qiong Fan, Xiaolan Hu, Suyun Jiang, Chunxiang Wang}}
\begin{document}
\maketitle
\begin{abstract}
For integers \(r\ge s\ge t\ge1\), we introduce the Ramsey number \(R(G;r,s,t)\), 
defined as the least \(N\) such that every assignment of an \(s\)-element subset of \([r]\) to each edge of \(K_N\) 
yields a copy of \(G\) with at least \(t\) colors common to all its edges. 
This parameter generalizes the classical set-coloring Ramsey number \(R(G;r,s)=R(G;r,s,1)\), 
whose study goes back to Erd\H{o}s, Hajnal, and Rado. We show that the parameter \(t\) cannot in general be eliminated by reparametrization: 
there are infinitely many quadruples \((G;r,s,t)\) for which \(R(G;r,s,t)\ne R(G;a,b)\) for all integers \(a\ge b\ge1\).

Write \(\alpha=\binom{r}{t}/\binom{s}{t}\). Using the inequality \(R(G;r,s,t)\le R(G;\binom rt,\binom st,1)\) 
together with the confirmed Erd\H{o}s–S\'os conjecture, we obtain \(R(T_n;r,s,t)\le\lfloor\alpha(n-2)\rfloor+2\) for every \(n\)-vertex tree \(T_n\). 
For stars, we give bounds that differ by one and determine the exact value under suitable divisibility conditions. 
For paths or even cycles \(G\) of order \(n\), we improve the upper bound to \(R(G;r,s,t)\le(\alpha-\gamma+o(1))n\) as \(n\to\infty\), 
for every fixed \(r>s\ge t\ge1\) and some positive constant \(\gamma=\gamma(r,s,t)\). 
We also show that the tree bound is asymptotically sharp for paths when \(r,n\to\infty\), 
and for even cycles in the iterated limit in which \(n\to\infty\) first and then \(r\to\infty\), with \(s,t\) fixed.

For cliques \(K_n\), we transfer known bounds of Conlon et al. for the ordinary set-coloring Ramsey numbers 
\(R(K_n;r,s)\) to the general intersection parameter \(t\), obtaining two-sided exponential bounds. 
For each fixed \(n\ge3\), we identify a nonnegative function \(\Phi_n\) with the following property. 
If \(s/r\to\sigma\) and \(t/r\to\tau\) as \(r\to\infty\), where \(0<\tau\le\sigma<1\), 
then \(R(K_n;r,s,t)=O(1)\) when \(\tau<\Phi_n(\sigma)\), whereas \(R(K_n;r,s,t)\ge\exp(cr)\) for some \(c>0\) when \(\tau>\Phi_n(\sigma)\).

\end{abstract}
\medskip
\noindent\textbf{Keywords.} Set-coloring Ramsey number; common intersection; tree; packing; path; even cycle; clique.
\par\smallskip
\noindent\textbf{MSC 2020.} 05C55; 05C35; 05B05.

\section{Introduction}
\label{sec:intr}

All graphs are finite and simple.
Write $[r]=\{1,\ldots,r\}$ and let $\binom{[r]}s$ denote its
family of $s$-subsets. An \emph{$(r,s)$-set-coloring} of a graph
$F$ is a map $\phi:E(F)\to\binom{[r]}s$.
The set-coloring Ramsey number $R(G;r,s)$ is the least $N$ such
that every such coloring of $K_N$ contains a copy of $G$ with a
color common to all its edges. This problem goes back to
Erd\H{o}s, Hajnal and Rado~\cite{EHR}. It includes ordinary
multicolor Ramsey theory when $s=1$, and the color-avoidance
problem studied by Erd\H{o}s and Szemer\'edi~\cite{ES} when $s=r-1$.

\begin{samepage}
We prescribe the size of the common intersection.
\begin{definition}
For integers $r\ge s\ge t\ge1$ and a graph $G$ with at least one edge, 
let $R(G;r,s,t)$ be the least $N$ such that every map
$\phi:E(K_N)\to\binom{[r]}s$ contains a \emph{$t$-monochromatic} copy $H\cong G$, 
meaning that
\(
 \left|\bigcap_{e\in E(H)}\phi(e)\right|\ge t.
\)
\end{definition}
\end{samepage}
The same $t$ colors must occur on every edge of the copy.
Writing $R(G; r)$ for the ordinary $r$-color Ramsey number, 
we have
\begin{equation}\label{eq:endpoints}
 R(G;r,s,1)=R(G;r,s),\ \ \ R(G;r,s,s)=R\left(G;\binom rs\right). 
\end{equation}
By default, $R(G;r,r,t)=|V(G)|$. 
We write $P_n$ and $C_n$ for the path and cycle on $n$ vertices, 
and $K_{1,n}$ for the star with $n$ leaves.

The parameter \(t\) cannot in general be eliminated by
changing the other two parameters. In
Section~\ref{sec:separation}, we show that the set of
values of \(R(P_3;r,s,t)\) is precisely the set of all
integers at least \(3\) except \(4\), whereas
\(R(P_3;a,b,1)\) is always odd. In
Section~\ref{sec:stars}, Example~\ref{star-gap} provides
a complementary obstruction for stars: for the fixed
triple \((r,s,t)=(4,3,2)\), no pair \(a\ge b\ge1\)
satisfies \( R(K_{1,n};4,3,2)=R(K_{1,n};a,b,1)\) 
simultaneously for all \(n\ge2\).

Replacing each edge label $S$ by the family $\binom{S}{t}$ of all
its $t$-subsets gives the shadow inequality
\begin{equation}\label{shadow1}
 R(G;r,s,t)\le R\left(G;\binom rt,\binom st,1\right).
\end{equation}
In Section~2, we prove that
$R(P_3;6,3,2)=6<7=R(P_3;15,3,1)$, showing that this inequality
can be strict. We also obtain lower bounds through a complementary
packing reduction, replacing ordinary colors with $s$-subsets of
$[r]$ whose pairwise intersections have size less than $t$.

Let $\alpha:=\binom rt/\binom st$. Combining the shadow reduction
\eqref{shadow1} with the recently confirmed Erd\H{o}s--S\'os
conjecture~\cite{RS} gives
\begin{equation}\label{eq:intro-tree-bound}
 R(T_n;r,s,t)\le\left\lfloor\alpha(n-2)\right\rfloor+2
\end{equation}
for every $n$-vertex tree $T_n$. We examine the sharpness of this
bound for stars and paths, and study related bounds for even cycles.

For stars, we relate $R(K_{1,n};r,s,t)$ to a local multiset
packing number. Once this packing number is known, the Ramsey
number is restricted to two consecutive values, with
factorization arguments settling several parity cases.
Estimates for the packing number give, for fixed $r,s,t$,
\[
 R(K_{1,n};r,s,t)=\alpha(n-1)+O_{r,s,t}(1)
 \qquad(n\to\infty),
\]
so the coefficient $\alpha$ in the tree bound is optimal.
R\"odl's packing theorem~\cite{Rodl}, with the case $t=s$
following directly, also gives
\[
 R(K_{1,n};r,s,t)=(1+o(1))\alpha(n-1)
\]
as $r\to\infty$ with $s,t$ fixed, uniformly for $n\ge2$.
When a Steiner system $S(t,s,r)$ exists, the local packing number equals $\alpha(n-1)$, 
and hence the Ramsey number is either $\alpha(n-1)+1$ or $\alpha(n-1)+2$.
Moreover, when $r-s\le t$ and $\binom{r-t}{s-t}\mid(n-1)$, 
an incidence-rigidity argument yields exact values by showing that 
every label class must be regular in any coloring of $K_{\alpha(n-1)+1}$ 
that avoids a $t$-monochromatic $K_{1,n}$.

For paths and even cycles, we obtain a strict improvement when the color parameters are fixed. 
For every fixed $r>s\ge t\ge1$, there is a constant
$\gamma=\gamma(r,s,t)>0$ such that, as $n\to\infty$,
\[
 R(P_n;r,s,t)\le(\alpha-\gamma+o(1))n,
 \qquad
 R(C_n;r,s,t)\le(\alpha-\gamma+o(1))n
 \quad(n\text{ even}).
\]
Thus, for paths with fixed color parameters, the upper bound
\eqref{eq:intro-tree-bound} has a linear gap from the true value.

When $s,t$ are fixed and both $r$ and $n$ tend to infinity,
the same upper bound is asymptotically sharp for paths:
\[
 R(P_n;r,s,t)=(1+o(1))\alpha n.
\]
For even cycles, we prove the corresponding conclusion in an
iterated limit: first $n\to\infty$ through even integers,
and then $r\to\infty$.

We also consider the regime in which $r,s,n\to\infty$,
$t$ is fixed, and $r/s\to\beta>1$.
By Theorem~\ref{thm:path-dense-asymptotic}, if
$\beta=m\ge2$ is an integer and $\log r=o(n)$, then
\[
 R(P_n;r,s,t)=(m^t+o(1))n=(1+o(1))\alpha n.
\]
In fact, this asymptotic formula holds with $P_n$ replaced by any $n$-vertex tree, 
so the tree bound is asymptotically sharp throughout this regime.
For paths, when $\beta^t$ is not an integer, we obtain a positive linear gap from $\alpha n$.
When $t=1$, we determine the leading coefficient for every limiting ratio $\beta>1$, 
without any restriction on the relative growth of $r$ and $n$.

Finally, for cliques, we derive general exponential upper and lower bounds for \(R(K_n;r,s,t)\) 
by transferring known bounds for the ordinary set-coloring Ramsey numbers.  
Combining the shadow upper bound \eqref{shadow1}, the elementary lower-bound comparison
in Proposition~\ref{prop:padding} and disjoint-witness reduction in Proposition~\ref{prop:witness},  
with the results of Conlon, Fox, He, Mubayi, Suk, and Verstra\"ete~\cite{CFHMSV},
we obtain two-sided bounds for fixed integers \(r>s\ge t\ge1\), and an improved upper bound in the regime \(\lfloor r/t\rfloor>r-s\).

We also investigate a different asymptotic regime in which the clique order \(n\) is fixed and the palette size \(r\) tends to infinity. 
For each fixed \(n\ge3\), we identify an explicit piecewise-linear function \(\Phi_n\) with the following property. 
If $s/r\to\sigma$ and $t/r\to\tau$ as $r\to\infty$, where $0<\tau\le\sigma<1$, 
then $R(K_n;r,s,t)$ remains bounded when $\tau<\Phi_n(\sigma)$ and grows 
at least exponentially in $r$ when $\tau>\Phi_n(\sigma)$.

Related questions for paths, stars, and other sparse graphs were considered in a preprint of He and Mao~\cite{HM}.
Buci\'c and Khamseh~\cite{BK} determined the exact Ramsey number for a path 
that avoids at least one color in an ordinary edge-coloring. 
Connections between set-coloring Ramsey numbers and error-correcting codes are developed in~\cite{CFHMSV,CFPZ}.
When $s=r-1$, assigning each edge its unique missing color identifies our problem with that of forcing a copy of $G$
using at most $r-t$ colors in an ordinary $r$-coloring.
Our results address both this color-avoidance case and general set-colorings.

\section{The intersection parameter cannot be eliminated}
\label{sec:separation}

In this section we show that the intersection parameter
$t$ cannot be eliminated. We prove that, for $P_3$,
there are infinitely many parameter triples whose Ramsey
numbers cannot be represented by any set-coloring Ramsey
number with one common color
(Corollary~\ref{cor:p3-separation}).
Moreover, the full set of values attained by
$R(P_3;r,s,t)$ consists of all integers at least $3$
except $4$ (Theorem~\ref{thm:p3-spectrum}).
In particular, the additional values obtained by allowing
$t\ge2$, beyond those already attained for $t=1$,
are exactly the even integers at least $6$.

\begin{lemma}\label{lem:p3-one}
For all integers $a\ge b\ge1$,
\begin{equation}\label{eq:p3-one}
 R(P_3;a,b,1)=2\left\lfloor\frac{a+b}{2b}\right\rfloor+1.
\end{equation}
In particular, $R(P_3;a,b,1)$ is always odd.
\end{lemma}
\begin{proof}
Recall that $K_N$ is class~1 when $N$ is even and class~2 when $N\ge3$ is odd. 
Thus, for $N\ge2$, its edge-chromatic number is
\[
 \chi'(K_N)=
 \begin{cases}N-1,&N\text{ even},\\N,&N\text{ odd}.\end{cases}
\]
An $(a,b)$-set-coloring of $K_N$ without a one-common-color
$P_3$ exists precisely when $a\ge b\chi'(K_N)$. 
For sufficiency, replace each color in a proper edge-coloring by a disjoint $b$-set. 
For necessity, every individual color must form a matching, so
\[
 b\binom N2\le a\left\lfloor\frac N2\right\rfloor,
\]
which gives the same condition. 
The largest $N$ with $a\ge b\chi'(K_N)$ is $2\lfloor(a+b)/(2b)\rfloor$, and \eqref{eq:p3-one} follows.
\end{proof}

The next theorem turns a classical design-theoretic object
--- coverings of two-edge paths by Hamilton cycles --- into exact values of $R(P_3;r,s,t)$.

\begin{theorem}\label{thm:p3-cover}
Let $q\ge3$ and $\lambda\ge1$ be integers, and suppose that $K_{q+1}$ admits a multiset $\mathcal H$ of Hamilton cycles
in which every two-edge path belongs to exactly $\lambda$ members of $\mathcal H$. Then
\begin{equation}\label{eq:p3-cover}
 R\!\left(P_3;\lambda\binom q2,\,\lambda(q-1),\,
 \lambda+1\right)=q+2.
\end{equation}
\end{theorem}
\begin{proof}
For the lower bound, regard the members of $\mathcal H$,
counted with multiplicity, as the colors. 
At a fixed vertex $v$ there are $\binom q2$ two-edge paths centered at $v$, 
and each Hamilton cycle contains exactly one of them, so the number of colors is $r=\lambda\binom q2$. 
Label each edge $uv$ of $K_{q+1}$ by the cycles containing it. 
There are $q-1$ two-edge paths centered at $u$ and containing $uv$;
each occurs in $\lambda$ cycles, and each cycle containing $uv$ accounts for exactly one of these paths. 
Hence every edge receives exactly $s=\lambda(q-1)$ colors. 
Two adjacent edges $uv$ and $vw$ have as common colors precisely the cycles containing the two-edge path $u$-$v$-$w$, 
of which there are exactly $\lambda$. 
Therefore this coloring of $K_{q+1}$ contains no $P_3$ with $\lambda+1$ common colors,
and $R(P_3;r,s,\lambda+1)\ge q+2$.

For the upper bound, consider any $(r,s)$-set-coloring of $K_{q+2}$ and a vertex with $m$ incident labels 
whose pairwise intersections have size at most $\lambda$. 
Let $d_i$ be the number of these labels containing color~$i$. 
Then $\sum_i d_i=ms$, and double counting the pairs of labels sharing color~$i$ gives, by Cauchy--Schwarz,
\[
 \frac{m^2s^2}{r}-ms
 \le\sum_{i=1}^{r}d_i(d_i-1)
 \le\lambda m(m-1),
\]
hence
\[
 m\left(\frac{s^2}{r}-\lambda\right)\le s-\lambda.
\]
Substituting $r=\lambda\binom q2$ and $s=\lambda(q-1)$, 
we have $s^2/r-\lambda=\lambda(q-2)/q$ and $s-\lambda= \lambda(q-2)$, so $m\le q$. 
But every vertex of $K_{q+2}$ has $q+1$ incident edges, 
so two of its incident labels share at least $\lambda+1$ colors. 
Thus $R(P_3;r,s,\lambda+1)\le q+2$, proving~\eqref{eq:p3-cover}.
\end{proof}

\begin{corollary}\label{cor:p3-separation}
For every $q=2^k$ with $k\ge2$,
\[
 R\!\left(P_3;\tbinom q2,\,q-1,\,2\right)=q+2.
\]
Since $q+2$ is even, Lemma~\ref{lem:p3-one} implies that none
of these values can be represented by a set-coloring Ramsey
number with one common color.
\end{corollary}
\begin{proof}
Nakamura, Kiyasu-Zen'iti and Ikeno~\cite{NKI80} constructed
an exact covering of the two-edge paths of $K_{q+1}$ by Hamilton cycles for every prime power $q$; 
apply Theorem~\ref{thm:p3-cover} with $\lambda=1$.
\end{proof}

The case $q=4$ is the example $R(P_3;6,3,2)=6$ from the introduction; the next cases are
\[
 R(P_3;28,7,2)=10,\qquad
 R(P_3;120,15,2)=18,\qquad
 R(P_3;496,31,2)=34.
\]

Taking $\lambda=2$ shows that, in fact, every even value occurs, already with $t=3$ fixed.

\begin{theorem}\label{thm:p3-spectrum}
The values attained by $R(P_3;r,s,t)$, as $r\ge s\ge t\ge1$ vary, 
are exactly the odd integers at least~$3$ and the even integers at least~$6$.
\end{theorem}
\begin{proof}
Every odd integer $2j+1\ge3$ equals $R(P_3;2j,1,1)$ by~\eqref{eq:p3-one}. 
For the even values, Kobayashi, Mutoh, Kiyasu-Zen'iti and Nakamura~\cite{KMKZN02} proved that every
complete graph of odd order has a multiset of Hamilton cycles in which each two-edge path occurs exactly twice. 
Taking $\lambda=2$ in Theorem~\ref{thm:p3-cover} gives
\[
 R\!\left(P_3;\,q(q-1),\,2(q-1),\,3\right)=q+2
 \qquad\text{for every even }q\ge4,
\]
realizing every even integer at least~$6$. 
Finally, the value~$4$ cannot occur for any $(r,s,t)$: 
if an $(r,s)$-set-coloring of $K_3$ avoids a $t$-common-color $P_3$, 
then its three labels intersect pairwise in fewer than $t$ colors, 
and assigning these three labels to the three perfect matchings of $K_4$ yields an avoiding coloring of $K_4$ as well.
\end{proof}

\begin{remark}\label{rem:shadow-gap}
The family in Corollary~\ref{cor:p3-separation} also shows that the shadow gap is multiplicatively unbounded, 
even for the fixed graph $P_3$ and fixed $t=2$. Put $r=\binom q2$ and $s=q-1$. Then
\[
 \frac{\binom r2}{\binom s2}=\frac{q(q+1)}4,
\]
which is an even integer for $q=2^k\ge8$, and
\eqref{eq:p3-one} gives
\[
 R\!\left(P_3;\tbinom r2,\tbinom s2,1\right)
 =\frac{q(q+1)}4+1,
 \qquad
 R(P_3;r,s,2)=q+2.
\]
The ratio of the shadow estimate to the true Ramsey number is asymptotic to $q/4$; the first two cases after
$R(P_3;6,3,2)=6<7=R(P_3;15,3,1)$ are $10<19$ and $18<69$.
\end{remark}

\section{Reductions and the tree bound}
\label{sec:reductions}

\begin{proposition}\label{prop:witness}
Let $\cW$ be a nonempty family of $t$-subsets of $[r]$, 
and suppose that every $s$-subset of $[r]$ contains at least $b\ge1$ members of $\cW$. 
Then, for every graph $G$,
\[
 R(G;r,s,t)\le R(G;|\cW|,b,1).
\]
In particular, the proposition yields \eqref{shadow1} and the following inequality.
\begin{align}
 R(G;r,s,t)&\le R(G;m,m-d,1)
 \quad\text{if }m=\lfloor r/t\rfloor>d=r-s.\label{eq:block}
\end{align}
\end{proposition}

\begin{proof}
Replace each label $S$ by any $b$ members of $\cW$ contained in $S$.
A common color in the resulting set-coloring is a common $t$-set in the original coloring. 
For~\eqref{shadow1}, take $\cW=\binom{[r]}t$. 
For~\eqref{eq:block}, take $m$ disjoint $t$-sets. 
The $d$ elements outside an $s$-set can meet at most $d$ of these blocks, 
so at least $m-d$ blocks lie inside it.
\end{proof}
We call \eqref{eq:block} the \emph{disjoint-witness reduction},
since its witness family consists of pairwise disjoint $t$-sets.

\begin{proposition}\label{prop:padding}
For every graph $G$,
\[
 R(G;r-t+1,s-t+1,1)\le R(G;r,s,t).
\]
\end{proposition}

\begin{proof}
To a set-coloring on a palette of size $r-t+1$, adjoin the same
$t-1$ new colors to every edge label. The intersection of the new
labels on a copy of $G$ has size at least $t$ only if the old labels
have a common color.
\end{proof}

When $s=r-1$, there is also a missing-color interpretation.
If $\rho_k^r(G)$ denotes
the least $N$ forcing a copy of $G$ using at most $k$ colors in every
ordinary $r$-coloring of $K_N$, then taking the unique missing color
of each edge gives
\begin{equation}\label{eq:missing}
 R(G;r,r-1,t)=\rho_{r-t}^r(G).
\end{equation}

Let $\nu(r,s,t)$ be the largest size of a family
$\cF\subseteq\binom{[r]}s$ whose distinct members intersect in fewer
than $t$ elements. Equivalently, $\cF$ is a $t$-packing. Then
\begin{equation}\label{eq:packing-count}
 \nu(r,s,t)\le\lfloor\alpha\rfloor.
\end{equation}

\begin{proposition}\label{prop:packing}
For every graph $G$,
\[
 R(G;r,s,t)\ge R(G;\nu(r,s,t)).
\]
\end{proposition}

\begin{proof}
Replace the colors in a Ramsey-critical ordinary coloring by the
members of a maximum packing. A common $t$-set can occur in only one
packing member, so a $t$-monochromatic copy would have been
monochromatic in the original coloring.
\end{proof}

Two standard packing results will be useful. R\"odl's theorem
\cite{Rodl} gives, for fixed $1\le t<s$,
\begin{equation}\label{eq:rodl}
 \nu(r,s,t)=(1-o(1))\alpha\qquad(r\to\infty).
\end{equation}
For $t=s$ equality holds in~\eqref{eq:packing-count}.
If a Steiner system $S(t,s,r)$ exists, then $\nu(r,s,t)=\alpha$
exactly. By the design existence theorem~\cite{Keevash,KeevashShort,GKLO}, for
fixed $t<s$ this holds for all sufficiently large $r$ satisfying
\begin{equation}\label{eq:design-div}
 \binom{s-i}{t-i}\ \bigm|\ \binom{r-i}{t-i}
 \qquad(0\le i<t).
\end{equation}
The asymptotic statement~\eqref{eq:rodl} requires no divisibility.

For a set-coloring $\phi$, let $H_T$ be the spanning graph with
\[
 E(H_T)=\{e:T\subseteq\phi(e)\}
 \qquad\left(T\in\binom{[r]}t\right).
\]
We use these support graphs throughout.

\begin{lemma}\label{lem:extremal}
If $\binom N2/\alpha>\ex(N,G)$, then $R(G;r,s,t)\le N$.
Consequently, if $\ex(N,G)\le cN/2$ for every $N$, then
\[
 R(G;r,s,t)\le\lfloor\alpha c\rfloor+2.
\]
\end{lemma}

\begin{proof}
Each edge belongs to $\binom st$ support graphs, so
$\sum_T e(H_T)=\binom st\binom N2$. Some support has at least
$\binom N2/\alpha$ edges. For the second assertion take
$N=\lfloor\alpha c\rfloor+2$.
\end{proof}

This is the usual extremal argument for set-colorings; see, for
example, He and Mao~\cite{HM}. The Erd\H{o}s--Gallai theorem
\cite{EG} gives $\ex(N,P_n)\le(n-2)N/2$. The confirmed Erd\H{o}s--S\'os
conjecture gives the same bound for every $n$-vertex tree $T_n$. (See Riordan and Scott~\cite{RS}.) Hence
\begin{equation}\label{eq:tree-bound}
 R(T_n;r,s,t)\le\lfloor\alpha(n-2)\rfloor+2.
\end{equation}
For paths this conclusion uses only the classical
Erd\H{o}s--Gallai theorem.

\section{Stars: local packings and rigidity}
\label{sec:stars}

For stars, the error in the tree bound is controlled by the
realizability of local packing data. For an integer \(h\ge0\), let
\(\nu(r,s,t;h)\) be the maximum size, counted with multiplicity, of
a multiset of \(s\)-subsets of \([r]\) in which each \(t\)-set belongs
to at most \(h\) members. Thus \(\nu(r,s,t;1)=\nu(r,s,t)\), and
incidence counting gives
\begin{equation}\label{eq:multi-count}
 \nu(r,s,t;h)\le\lfloor\alpha h\rfloor.
\end{equation}

\subsection{Local packing bounds and designs}
\label{subsec:star-local}

\begin{theorem}\label{thm:star-packing}
For \(n\ge2\),
\[
 \nu(r,s,t;n-1)+1\le R(K_{1,n};r,s,t)
 \le\nu(r,s,t;n-1)+2.
\]
The upper value holds if \(\nu(r,s,t;n-1)\) is odd, or if a maximum
multiset has even multiplicity for every label.
\end{theorem}
\begin{proof}
Put \(D=\nu(r,s,t;n-1)\). In a coloring avoiding the desired star,
the incident labels at each vertex form a feasible multiset, so
its degree is at most \(D\).
If \(D\) is odd, assign a maximum multiset to the \(D\) factors of a
\(1\)-factorization of \(K_{D+1}\). If \(D\) is even, delete one label
and use a \(1\)-factorization of \(K_D\). These constructions give the
lower bounds. If all multiplicities are even, decompose \(K_{D+1}\)
into \(2\)-factors and assign half the multiplicity of each label to
that many factors. At every vertex the incident multiset is then
the prescribed maximum multiset.
\end{proof}

The general bound in Theorem~\ref{thm:star-packing} is nearly tight.
To see this, write $n-1=q\binom{r-t}{s-t}+h$, where 
$0\le h<\binom{r-t}{s-t}$. Adjoining \(q\) copies of every \(s\)-set 
to a feasible residual multiset gives
\begin{equation}\label{eq:residual}
 \nu(r,s,t;n-1)\ge q\binom rs+\nu(r,s,t;h).
\end{equation}
In particular,
\begin{equation}\label{eq:star-asymptotic}
 \binom rs\left\lfloor\frac{n-1}{\binom{r-t}{s-t}}\right\rfloor+1
 \le R(K_{1,n};r,s,t)
 \le\lfloor\alpha(n-1)\rfloor+2.
\end{equation}
Hence, for fixed \(r,s,t\),
\begin{equation}\label{eq:star-slope}
 R(K_{1,n};r,s,t)=\alpha(n-1)+O_{r,s,t}(1).
\end{equation}
The shadow-extremal upper bound is therefore asymptotically sharp
for stars as \(n\to\infty\). 

When the palette is large compared to the prescribed intersection,
the local packing number can be computed exactly via designs.

\begin{corollary}\label{cor:star-design}
For \(n\ge2\),
\[
 (n-1)\nu(r,s,t)+1\le R(K_{1,n};r,s,t)
 \le\lfloor\alpha(n-1)\rfloor+2.
\]
For fixed \(s\ge t\ge1\), as \(r\to\infty\),
\[
 R(K_{1,n};r,s,t)=(1+o(1))\alpha(n-1)
\]
uniformly for integers \(n\ge2\).
If a Steiner system \(S(t,s,r)\) exists, then
\[
 R(K_{1,n};r,s,t)\in
 \{\alpha(n-1)+1,\,\alpha(n-1)+2\}.
\]
The upper value holds, and the shadow bound is exact, if
\(\alpha(n-1)\) is odd or \(n\) is odd.
\end{corollary}
\begin{proof}
Repeat each block of a maximum packing \(n-1\) times to obtain
\(\nu(r,s,t;n-1)\ge(n-1)\nu(r,s,t)\), and apply
Theorem~\ref{thm:star-packing}. 
For \(t<s\), R\"odl's packing theorem~\cite{Rodl} gives
\(\nu(r,s,t)/\alpha\to1\); for \(t=s\), we have
\(\nu(r,s,s)=\alpha\) directly. Thus
\(\nu(r,s,t)/\alpha\to1\) for all fixed \(s\ge t\).
Moreover, for every \(n\ge2\),
\[
\frac{\nu(r,s,t)}{\alpha}
\le
\frac{R(K_{1,n};r,s,t)}{\alpha(n-1)}
\le
1+\frac{2}{\alpha(n-1)}
\le
1+\frac{2}{\alpha}.
\]
Since \(\nu(r,s,t)/\alpha\to1\) and \(\alpha\to\infty\),
the convergence is uniform over all integers \(n\ge2\).
For a Steiner system, \(\nu(r,s,t)=\alpha\), so the local packing
bound is attained. If \(n\) is odd, the repeated block multiplicities
are even. Theorem~\ref{thm:star-packing} therefore gives
\(R(K_{1,n};r,s,t)=\alpha(n-1)+2\) in both cases,
which also equals the extremal upper bound for the shadow Ramsey
number.
\end{proof}

\subsection{Incidence rigidity and exact values}
\label{subsec:star-rigidity}

\begin{theorem}\label{thm:rigid-star}
Suppose $r-s\le t$ and \(n-1=q\binom{r-t}{s-t}\) for an integer
\(q\ge1\). Then
\[
 R(K_{1,n};r,s,t)=
 \begin{cases}
 \alpha(n-1)+1,&\text{if $\binom rs$ is even and $q$ is odd},\\
 \alpha(n-1)+2,&\text{otherwise}.
 \end{cases}
\]
Moreover, in any coloring of \(K_{\alpha(n-1)+1}\) avoiding the
specified star, every \(s\)-set labels a \(q\)-regular spanning subgraph.
\end{theorem}

\begin{proof}
Write \(M=\binom rs\) and \(C=\binom{r-t}{s-t}\), so \(n-1=Cq\) and \(\alpha(n-1)=Mq\).
We split the proof into two claims.

\smallskip
\noindent\textit{Claim 1.} \(\alpha(n-1)+1\le R(K_{1,n};r,s,t)\le\alpha(n-1)+2\).

Taking each $s$-set $q$ times gives a feasible multiset
of size $Mq$. The counting bound shows that this multiset
is maximum. Hence Theorem~\ref{thm:star-packing} gives
\[
 Mq+1\le R(K_{1,n};r,s,t)\le Mq+2.
\]
Moreover, if $Mq$ is odd or $q$ is even, the same theorem
gives $R(K_{1,n};r,s,t)=Mq+2$.

\smallskip
\noindent\textit{Claim 2.} 
If \(r-s\le t\), let $x_S$ for \(S\in\binom{[r]}s\) be real numbers satisfying
\[
 \sum_{S\supseteq T}x_S=q\binom{r-t}{s-t}
 \qquad\text{for every }T\in\binom{[r]}t.
\]
Then \(x_S=q\) for every \(S\). 

To see this, for \(J\subseteq[r]\) write \(d(J)=\sum_{S\supseteq J}x_S\).
If \(|J|=j\le t\), sum the hypotheses over the \(t\)-sets containing
\(J\) to obtain
\[
 \binom{s-j}{t-j}d(J)
 =q\binom{r-j}{t-j}\binom{r-t}{s-t},
 \qquad d(J)=q\binom{r-j}{s-j}.
\]
Fix \(D\in\binom{[r]}{r-s}\). Since \(|D|=r-s\le t\),
inclusion--exclusion gives
\[
 x_{[r]\setminus D}
 =\sum_{J\subseteq D}(-1)^{|J|}d(J)
 =q\sum_{j=0}^{r-s}(-1)^j\binom{r-s}{j}\binom{r-j}{s-j}=q.
\]
The last sum counts the \(s\)-sets disjoint from \(D\), of which there
is exactly one (a binomial coefficient with negative lower index is zero).

\smallskip
Returning to the proof of the theorem. 
Suppose that a coloring of $K_{Mq+1}$ avoids the specified star.
For every vertex $v$ and every $T\in\binom{[r]}t$,
we have $d_{H_T}(v)\le n-1=Cq$. On the other hand,
\[
 \sum_{T\in\binom{[r]}t}d_{H_T}(v)
 =\binom st Mq
 =\binom rt Cq.
\]
Thus $d_{H_T}(v)=Cq$ for every $v,T$.
Applying Claim 2 to the incident multiplicities $x_S(v)$
shows that every label class is a $q$-regular spanning subgraph,
proving the structural assertion.

If $M$ is even and $q$ is odd, then $Mq+1$ is odd.
The handshaking lemma therefore rules out such a coloring.
Consequently $R(K_{1,n};r,s,t)\le Mq+1$,
and Claim 1 gives equality.
Otherwise, either $Mq$ is odd or $q$ is even,
so Claim 1 gives $R(K_{1,n};r,s,t)=Mq+2$.
\end{proof}

Theorem~\ref{thm:rigid-star} yields an infinite family of
strict shadow gaps.
\begin{example}\label{star-gap}
If \(n\equiv3\pmod4\), then
\begin{equation}\label{eq:infinite-shadow-gap}
 R(K_{1,n};4,3,2)=2n-1
 <R(K_{1,n};6,3,1)=2n.
\end{equation}
The first equality follows by taking \(q=(n-1)/2\) in
Theorem~\ref{thm:rigid-star}. For the second, split a
\(2\)-factorization of \(K_{2n-1}\) into two equal
collections of \(2\)-factors. This gives two
\((n-1)\)-regular spanning graphs, which we label by
two disjoint \(3\)-subsets of \([6]\).
Every \(n\)-edge star contains an edge from each graph,
so its edge labels have empty common intersection.
Thus \(R(K_{1,n};6,3,1)\ge 2n\). The reverse inequality
follows from \eqref{eq:star-asymptotic}.

In fact, no fixed pair of parameters for the ordinary
set-coloring problem can represent
\(R(K_{1,n};4,3,2)\) for all stars. More precisely,
there are no integers \(a\ge b\ge1\) such that
\[
 R(K_{1,n};4,3,2)=R(K_{1,n};a,b,1)
 \qquad\text{for every }n\ge2.
\]
Indeed, the asymptotic formula~\eqref{eq:star-slope}
forces \(a/b=2\), so \(a=2b\).
Since a Steiner system \(S(1,b,2b)\) exists,
Corollary~\ref{cor:star-design} gives
\[
 R(K_{1,n};2b,b,1)=2n
\]
for odd \(n\). But if \(n\equiv3\pmod4\), then
\eqref{eq:infinite-shadow-gap} gives
\(R(K_{1,n};4,3,2)=2n-1\ne2n\), contradicting the assumed
equality for every \(n\ge2\).
\end{example}

\section{Paths and even cycles}
\label{sec:paths}

In the study of ordinary multicolor path Ramsey numbers, 
it is conjectured that the lower bound is the truth, i.e., \(R(P_n; r) = (r - 1 + o(1))n\). 
The trivial upper bound \(R(P_n; r) \le rn\) was only recently broken, 
with the best current bound due to Knierim and Su~\cite{KS}, 
who proved \(R(P_n; r) \le (r - 1/2 + o(1))n\). 
In close analogy to this classical breakthrough, 
our main objective in this section is to break the trivial shadow upper bound for set-coloring Ramsey numbers. 
While the shadow reduction merely yields the coefficient \(\alpha = \binom{r}{t}/\binom{s}{t}\), 
we establish a strict improvement by showing \(R(P_n;r,s,t) \le (\alpha - \gamma + o(1))n\). 

\subsection{Upper bounds}\label{subsec:upper-bounds}

For \(x>1\) define
\begin{equation}\label{eq:beta}
 b(x)=\frac{x\lfloor x\rfloor-\sqrt{x\lfloor x\rfloor(\lfloor x\rfloor+1-x)}}{x-1}.
\end{equation}
Thus $b(x)=x$ at integers, and $\lfloor x\rfloor<b(x)<x$ otherwise. 
The function arises from the exact path extremal number, rather
than its linear relaxation. If $N=a(n-1)+u$, $0\le u<n-1$, the
Faudree--Schelp theorem~\cite{FS} gives
\begin{equation}\label{eq:exact-path-ex}
 E_n(N):=\ex(N,P_n)=a\binom{n-1}{2}+\binom u2.
\end{equation}
Applying this sharper extremal estimate gives the finite bound
\begin{equation}\label{eq:finite-path-upper}
 R(P_n;r,s,t)\le
 \min\left\{N\ge n:\binom N2>\alpha E_n(N)\right\}.
\end{equation}
This follows by summing~\eqref{eq:exact-path-ex} over the
$\binom rt$ supports.

If $N/n\to z$, then \(
 \frac{2E_n(N)}{n^2}\longrightarrow
 \lfloor z\rfloor+(z-\lfloor z\rfloor)^2 \). 
Put $j=\lfloor\alpha\rfloor$. On $j\le z\le j+1$,
the equation \( j+(z-j)^2=\frac{z^2}{\alpha} \) 
has smaller root $z=b(\alpha)$.
For $z>b(\alpha)$ sufficiently close to $b(\alpha)$,
we have \( z^2>\alpha\bigl(j+(z-j)^2\bigr) \). 
The finite bound~\eqref{eq:finite-path-upper} therefore gives
$R(P_n;r,s,t)\le(b(\alpha)+o(1))n$.

We next prove that the counting coefficient $\alpha$ can always be
improved, including when $\alpha$ is an integer. The improvement
below is explicit but small; its purpose is to rule out asymptotic
saturation of all supports simultaneously. Put
\begin{equation}\label{eq:gap-constant}
 \gamma=
 \frac{\alpha-1}
 {4096\binom st(\alpha+1)^2\lceil2(\alpha+1)\rceil^{\binom rt}},
 \qquad c_* =\min\{b(\alpha),\alpha-\gamma\}.
\end{equation}

\begin{theorem}\label{thm:long-path-cycle}
Fix $r>s\ge t\ge1$. As $n\to\infty$,
\[
 R(P_n;r,s,t)\le(c_*+o(1))n,
 \qquad R(C_n;r,s,t)\le(c_*+o(1))n\quad(n\text{ even}).
\]
The error terms may depend on $r,s,t$. In particular, the upper
coefficient is strictly smaller than $\alpha$ for every fixed
nontrivial palette.
\end{theorem}

The proof also applies to unrestricted shadow set-colorings. In the ordinary-color
case, stronger estimates are known: for $k\ge4$, Knierim and
Su~\cite{KS} give coefficient $k-\tfrac12$. Here the loss estimate also applies to overlapping supports.

A matching is \emph{connected} if its edges lie in one connected
component; its size is its number of edges. We use the
Erd\H{o}s--Gallai matching bound~\cite{EG}: if a graph on $z$ vertices
has no matching of size $h/2$, where $h\ge4$ is even, then
\begin{equation}\label{eq:matching-ex}
 e(H)\le\max\left\{\binom{h-1}{2},
 \binom{h/2-1}{2}+(h/2-1)(z-h/2+1)\right\}.
\end{equation}
We will apply it only to components with $z\ge h$.

\begin{lemma}\label{lem:support-loss}
Fix $r>s\ge t\ge1$, and let $\gamma$ be as in
\eqref{eq:gap-constant}. There exist $\eta>0$ and $h_0$ such that,
for every even $h\ge h_0$, every $(r,s)$-set-coloring of a graph $G$
on $v\ge(\alpha-\gamma)h$ vertices with at least
$(1-\eta)\binom v2$ edges has a connected matching of size $h/2$
in some $t$-support.
\end{lemma}
\begin{proof}
For this proof put
\[
 Q=\binom rt,\qquad P=\binom st,\qquad
 A=\alpha+1,\qquad B=\lceil2A\rceil^Q.
\]
Let $H_1,\ldots,H_Q$ be the supports and suppose that none has
such a matching. In particular, none contains $P_h$. Define
\[
 D_i=\frac{hv}{2}-e(H_i),\qquad D=\sum_{i=1}^Q D_i.
\]
For a component $U$ of $H_i$ write
$D_i(U)=h|U|/2-e(H_i[U])$. Erd\H{o}s--Gallai for paths gives
$D_i(U)\ge0$, and $D_i=\sum_U D_i(U)$.
As every edge of $G$ belongs to exactly $P$ supports,
\begin{equation}\label{eq:total-loss}
 D=\frac{Qhv}{2}-Pe(G).
\end{equation}
We first show that $D$ cannot be too small.

Set
\[
 \varepsilon=\frac{Q-P}{64QB},\qquad
 d_* =\frac{\varepsilon}{16A^2}.
\]
Notice that $\gamma=\alpha d_*/(4P)$ and
$0<\varepsilon<1/4$. For the moment assume
\begin{equation}\label{eq:loss-contradiction}
 h\le v\le Ah,\qquad v\ge8B,\qquad D\le d_*v^2.
\end{equation}
A component with $z\le h/2$ vertices has
\[
 D_i(U)\ge hz/2-\binom z2\ge hz/4.
\]
Thus the total number of vertices in these small components,
summed over all supports, is at most $4D/h$.
For a component with $z>(1+\varepsilon)h$,
\eqref{eq:matching-ex} implies
\[
 e(H_i[U])\le\max\{h^2/2,hz/2-h^2/8\}.
\]
Consequently
\[
 D_i(U)\ge\min\{h(z-h)/2,h^2/8\}
 \ge\varepsilon h^2/2.
\]
There are at most $2D/(\varepsilon h^2)$ such component
occurrences across all supports. Since each has at most $v\le Ah$
vertices, their total order is at most $2AD/(\varepsilon h)$.

Delete every vertex which lies in a small component or a component
of order greater than $(1+\varepsilon)h$ in at least one support.
Let $W$ be the remaining set. The union bound gives
\begin{align*}
 |V(G)\setminus W|
 &\le\left(4+\frac{2A}{\varepsilon}\right)\frac D h
 \le\frac{4AD}{\varepsilon h}
 \le\frac{4A^2d_*}{\varepsilon}v=\frac v4.
\end{align*}
In particular, $|W|\ge v/2$. Every vertex in $W$ lies, in each
$H_i$, in an original component of order between $h/2$ and
$(1+\varepsilon)h$. There are at most $2v/h\le2A$ such
components in each support.

Assign to each vertex of $W$ the tuple of its $Q$ component
indices. This partitions $W$ into at most $B$ cells. Two vertices
in the same cell lie in the same component of every support.
Nevertheless, their pair is an edge of at most $P$ supports,
even if it is absent from $G$. Hence it is missing from at least
$Q-P$ supports. If the cell sizes are $w_1,\ldots,w_b$, $b\le B$,
Cauchy--Schwarz gives
\begin{equation}\label{eq:cell-deficit}
 \sum_{j=1}^b\binom{w_j}{2}
 \ge\frac12\left(\frac{|W|^2}{B}-|W|\right)
 \ge\frac{v^2}{16B},
\end{equation}
where the last inequality uses $|W|\ge v/2$ and $v\ge8B$.
Thus the total number of missing pairs inside these cells, counted
with their missing support, is at least $(Q-P)v^2/(16B)$.

There is also an upper bound for the same quantity. For every
retained original component $U$ of order $z\le(1+\varepsilon)h$,
\[
 \binom z2-e(H_i[U])
 =D_i(U)+\frac{z(z-1-h)}2
 \le D_i(U)+\frac{\varepsilon hz}{2}.
\]
Summing over the retained components and then over supports shows
that their total number of missing pairs is at most
\[
 D+\frac{Q\varepsilon hv}{2}
 \le\left(d_*+\frac{Q\varepsilon}{2}\right)v^2
 <\frac{Q-P}{16B}v^2.
\]
Indeed, $Q\varepsilon/2=(Q-P)/(128B)$ and
$d_*\le(Q-P)/(1024B)$. This contradicts~\eqref{eq:cell-deficit}.
We have proved that the first two conditions
in~\eqref{eq:loss-contradiction} force $D>d_*v^2$.

We now take $\eta=d_*/(2P)$ and choose $h_0$ sufficiently large.
The path extremal bound and the density assumption give
\[
 P(1-\eta)\binom v2\le\sum_i e(H_i)
 \le Q(h-2)v/2.
\]
Since $\eta<1/(2A)$, this implies $v\le Ah$ for sufficiently
large $h$. Also $v\ge(\alpha-\gamma)h>h$; increase $h_0$ so that
$v\ge8B$ and $P/(2v)\le d_*/4$.
Finally, by~\eqref{eq:total-loss},
\begin{align*}
 \frac D{v^2}
 &\le\frac12\left(\frac{Q}{\alpha-\gamma}
       -P(1-\eta)(1-1/v)\right)\\
 &\le\frac{P\gamma}{\alpha}+\frac{P\eta}{2}
       +\frac{P}{2v}
 \le\frac{3d_*}{4}.
\end{align*}
Here we used $\gamma<\alpha/2$ and
$(1-\gamma/\alpha)^{-1}\le1+2\gamma/\alpha$.
All the conditions in~\eqref{eq:loss-contradiction} now hold,
which is the required contradiction.
\end{proof}

The obstruction in the proof is global. Near equality in the edge
count would make most vertices belong to almost complete components
of bounded order in every support. The common refinement of these
component partitions then supplies many pairs which would have to
belong to nearly all supports, contrary to their prescribed
multiplicity $\binom st<\binom rt$.

To obtain cycles of a prescribed length, we use the
connected-matching embedding lemma: a connected
matching covering $h$ vertices of a bounded reduced graph, with
all its edges and connecting edges represented by regular pairs of
positive density and equal cluster size $L$, yields every even
cycle of length between a constant depending on the reduced graph
and $(1-o(1))hL$. This is the standard even-cycle form of
\L uczak's method; see~\cite{LSS,KS}. The regularity error tends to
zero before $L$ tends to infinity. Connecting paths use only a
bounded number of vertices per cluster, and the matching pairs
supply the remaining even lengths.

\begin{lemma}\label{lem:set-transfer}
Fix $r,s,t$ and $c\ge1$. Suppose that for every $\xi>0$ there are
$\eta>0$ and $h_0$ such that every $(r,s)$-set-coloring of a graph
on $v\ge(c+\xi)h$ vertices with at least
$(1-\eta)\binom v2$ edges has a connected matching of size $h/2$
in a $t$-support, for all even $h\ge h_0$.
Then $R(C_n;r,s,t)\le(c+o(1))n$ for even $n$.
\end{lemma}

\begin{proof}
Fix $\zeta>0$, let $N=\lceil(c+\zeta)n\rceil$, and write
$M=\binom rs$ for this proof.
Apply the simultaneous regularity lemma to the $M$ individual
label graphs of a coloring of $K_N$. Choose the regularity
parameter sufficiently small in terms of $r,s,t,\zeta,\eta$,
and choose the minimum number of clusters sufficiently large in
terms of $h_0$ and $\zeta$. We obtain equal clusters of size $L$,
an exceptional set of $o(N)$ vertices, and a reduced graph $F$
on $v$ clusters with $e(F)\ge(1-\eta)\binom v2$. Its edges are
pairs regular in every label graph. The number $v$ is bounded
independently of $n$.

On each edge of $F$, choose a label whose density across that
pair is at least $1/M$. This is an $(r,s)$-set-coloring of $F$.
Choose the largest even integer $h\le v/(c+\zeta/3)$.
The choice of the minimum number of clusters guarantees both
$h\ge h_0$ and \(  hL>\left(1+\rho\right)n \) 
for a constant $\rho>0$ depending only on $c,\zeta$; the loss of
two clusters and the exceptional set are absorbed by the gap
between $c+\zeta$ and $c+\zeta/3$.
The assumed reduced-graph statement supplies a matching of $h/2$
edges in a connected $t$-support, say with common set $T$.

For every edge of this reduced support, the original $T$-support
has density at least $1/M$. It is also regular: it is a union of
at most $M$ individual label graphs, so its regularity error is
at most $M$ times the original error. Thus all the matching and
connecting pairs satisfy the hypotheses of the connected-matching
embedding lemma in the \emph{same} graph $H_T$. Taking the error
small enough relative to $\rho$ and then $n$ large enough embeds
$C_n$ in $H_T$. All its edges contain $T$, as required.
\end{proof}

\begin{proof}[\textbf{Proof of Theorem~\ref{thm:long-path-cycle}}]
Lemmas~\ref{lem:support-loss} and~\ref{lem:set-transfer} give coefficient $\alpha-\gamma$.
For coefficient $b(\alpha)$, suppose that no support in a reduced
graph contains a connected matching of size $h/2$. Then every
support is $P_h$-free, and \(  e(F)\le\alpha E_h(v)\). 
For $v/h\to z$, this implies, up to a factor tending to one,
\[
 z^2\le\alpha\bigl(\lfloor z\rfloor+
                         (z-\lfloor z\rfloor)^2\bigr).
\]
The ratio of the expression in parentheses to $z^2$ is continuous
and decreasing for $z\ge1$ (its derivative on $(j,j+1)$ is
$2j(z-j-1)/z^3$). It follows that for each $\xi>0$ this inequality
is impossible when $v\ge(b(\alpha)+\xi)h$, the density error is
sufficiently small, and $h$ is sufficiently large. Equivalently,
one may first choose an induced subgraph on
$\lceil(b(\alpha)+\xi/2)h\rceil$ vertices with at least the same
edge density, by averaging. Lemma~\ref{lem:set-transfer} now gives
the cycle bound. For paths, use $P_n\subset C_n$ when $n$ is even
and $P_n\subset C_{n+1}$ when $n$ is odd.
\end{proof}

For dense palettes we can also use the disjoint-witness reduction.
The missing-color identity~\eqref{eq:missing} and the theorem of
Buci\'c and Khamseh~\cite[Theorem~1.1]{BK} give
\begin{equation}\label{eq:bk}
 R(P_n;k,k-1,1)=n+\left\lfloor\frac{n-2}{2^k-2}\right\rfloor
 \qquad(k,n\ge2).
\end{equation}
Combining this exact result with the two palette reductions yields
the following explicit comparison.

\begin{theorem}\label{thm:dense-path}
Let $n\ge2$, $1\le t\le r-1$, and $m=\lfloor r/t\rfloor$.
Then
\[
 R(P_n;r,r-1,t)\ge
 n+\left\lfloor\frac{n-2}{2^{r-t+1}-2}\right\rfloor.
\]
If $m\ge2$, then also
\[
 R(P_n;r,r-1,t)\le
 n+\left\lfloor\frac{n-2}{2^m-2}\right\rfloor.
\]
In particular, if $n<2^m$, then $R(P_n;r,r-1,t)=n$.
\end{theorem}

\begin{proof}
The padding reduction (Proposition~\ref{prop:padding}) gives the lower comparison with
$R(P_n;r-t+1,r-t,1)$. The disjoint-witness reduction gives the
upper comparison with $R(P_n;m,m-1,1)$. Substitute~\eqref{eq:bk}. 
For the last assertion the upper floor is zero, and every Ramsey
number for $P_n$ is at least $n$.
\end{proof}

For example, $R(P_n;12,11,2)=n$ for $2\le n\le63$.
The general extremal estimate would give only
$\lfloor\tfrac65(n-2)\rfloor+2$. The theorem illustrates why a
small witness family can be more effective than the full shadow.

\subsection{Packing lower bounds and asymptotics}
\label{subsec:packing-asymptotics}

We next turn to lower bounds and asymptotic consequences of the
packing reduction.  The main point is that the packing number
$\nu(r,s,t)$ is asymptotically equal to the counting parameter
\( \alpha \) 
when the palette grows with $s$ and $t$ fixed.  Combined with the
upper bounds above, this yields asymptotic equivalence for paths
and, in an iterated sense, for even cycles.

\begin{theorem}
\label{thm:packing-asymptotic}
Fix integers $s\ge t\ge1$.  Then the following hold.

\begin{enumerate}
\item[\textup{(i)}] If $r,n\to\infty$, then
\[
R(P_n;r,s,t)
=
(1+o(1))\,
\frac{\binom rt}{\binom st}\,n.
\]
Moreover,
\[
\frac{
R(P_n;r,s,t)}
{
R\!\left(P_n;\binom rt,\binom st,1\right)
}
\longrightarrow1.
\]

\item[\textup{(ii)}] Fix $s\ge t\ge1$. As $r\to\infty$,
\[
 \liminf_{\substack{n\to\infty\\n\text{ even}}}
 \frac{R(C_n;r,s,t)}{\alpha n}\longrightarrow1,
 \   \limsup_{\substack{n\to\infty\\n\text{ even}}}
 \frac{R(C_n;r,s,t)}{\alpha n}\longrightarrow1.
\]
The same iterated asymptotic equivalence holds with the shadow
Ramsey number $R(C_n;\binom rt,\binom st,1)$. 
\end{enumerate}
\end{theorem}

\begin{proof}
Put \( k=\nu(r,s,t)\). 
By R\"odl's packing theorem when \(t<s\), and directly
when \(t=s\), we have \(k/\alpha\to1\) as \(r\to\infty\). 

(i) For paths,~\eqref{eq:tree-bound} gives  
\begin{equation}\label{eq:path-packing-upper}
R(P_n;r,s,t)\le \lfloor\alpha(n-2)\rfloor+2.
\end{equation}
For the lower bound, we recall the standard multicolor path construction
\cite{FS,SYXL}, giving its short description to fix rounding.
Put \( a=\left\lfloor\frac n2\right\rfloor-1\). 
Assume first that $k\ge2$ and $a\ge1$.  On
$\mathbb Z_{2k-2}$, for each
$i=0,\ldots,k-2$, assign color $i+1$ to the edges from $i$ to
\( i+1,\ldots,i+k-2 \) 
and from $i+k-1$ to \(i+k,\ldots,i+2k-3\), 
with indices taken modulo $2k-2$.  These form two vertex-disjoint
stars.  Assign the remaining antipodal matching color $k$.

Replace every vertex by a set of $a$ vertices, retain the colors
between parts, and assign all edges inside a part color $k$.
For colors $1,\ldots,k-1$, every component is bipartite with one
part of size $a$, and hence contains no path on $n$ vertices.
The color-$k$ components have at most $2a<n$ vertices.
Thus this gives a coloring of \( K_{2(k-1)a}\) 
with no monochromatic $P_n$.  By Proposition~\ref{prop:packing}, 
\begin{equation}\label{eq:path-packing-lower}
R(P_n;r,s,t)
\ge
\max\left\{
n,\,
2(k-1)\left(
\left\lfloor\frac n2\right\rfloor-1
\right)+1
\right\}.
\end{equation}
The exceptional cases give only the trivial bound $R(P_n;r,s,t)\ge n$.

Dividing~\eqref{eq:path-packing-lower} by
$\alpha n$ gives
\[
\frac{R(P_n;r,s,t)}{\alpha n}
\ge
\frac{k-1}{\alpha}\,
\frac{2(\lfloor n/2\rfloor-1)}{n}
+o(1),
\]
and hence the right-hand side tends to $1$ whenever
$r,n\to\infty$.  On the other hand,~\eqref{eq:path-packing-upper}
gives
\[
\frac{R(P_n;r,s,t)}{\alpha n}
\le
1+o(1).
\]

For the second assertion in~\textup{(i)}, 
the shadow Ramsey number lies between the original
one and $\lfloor\alpha(n-2)\rfloor+2$, by
Lemma~\ref{lem:extremal} applied to the shadow palette.

(ii) For even $n\ge4$ and $k \ge2$, we have
\begin{equation}\label{eq:cycle-lower}
R(C_n;r,s,t)\ge
\max\left\{
\frac{3n}{2}-1,\,
(k-1)(n-2)+2
\right\}.
\end{equation}
Indeed, the first construction consists of two red cliques of
orders $n-1$ and $n/2-1$, with all cross-edges colored blue.
Neither color contains a copy of $C_n$.  For the second,
use the $k$-color construction in the proof of~\eqref{eq:path-packing-lower} 
and add one vertex joined to every old vertex in color $k$. 
A cycle in color $k$ lies in one old component together with the new vertex, 
so has at most $n-1$ vertices. 
In every other color, each component is bipartite with one part of size $n/2-1$. 
The packing reduction completes the proof. If $k=1$, retain the lower bound $n$.

Using~\eqref{eq:cycle-lower}, we obtain 
\[
\liminf_{\substack{n\to\infty\\ n\ {\rm even}}}
\frac{R(C_n;r,s,t)}{\alpha n}
\ge
\frac{k-1}{\alpha}.
\]
Since 
\(\frac{b(\alpha)}{\alpha}\longrightarrow1\) 
and \(\frac{\gamma}{\alpha}\longrightarrow0\) as \(r\to\infty\), 
by Theorem~\ref{thm:long-path-cycle}, we have 
\begin{equation}\label{eq:cycle-limsup}
\lim_{r\to\infty}
\limsup_{\substack{n\to\infty\\ n\ {\rm even}}}
\frac{R(C_n;r,s,t)}{\alpha n}
\le1.
\end{equation}

For the shadow number use its lower bound by the original number
and the same upper argument for the shadow palette.
The palette dependence of the cycle embedding threshold dictates
the order of the two limits.
\end{proof}

The two-clique construction also improves the path lower bound
to $n+\lfloor n/2\rfloor-1$ whenever $k\ge2$.
For fixed parameters admitting a Steiner system $S(t,s,r)$,
$k=\alpha$; hence for long paths and even cycles the general
interval of coefficients is now
\[
 [\alpha-1,\alpha-\gamma],
\]
with the two-color lower bound used when it is larger. 
Determining the optimal coefficient within this interval includes
longstanding problems in the ordinary multicolor subcase.

For example, $(r,s,t)=(5,3,2)$ gives
\[
 c_*\le b(10/3)=\frac{30-2\sqrt{15}}7<3.1792,
\]
a substantial improvement over the counting coefficient
$\alpha=10/3$.

For $s=r-1$ and $r>2t$, the exact-extremal coefficient is
\[
 b(\alpha)=\frac{2}{1+\sqrt{1-2t/r}}
 =1+\frac{t}{2r}+O((t/r)^2).
\]
The even-cycle estimate takes $n\to\infty$ with the palette fixed;
the final expansion describes this coefficient as $t/r\to0$.
The disjoint-witness path bound can be stronger still.

\subsection{Dense-palette asymptotics for paths}
\label{subsec:path-dense-asymptotics}

The upper bound in the following theorem follows from the finite
extremal inequality~\eqref{eq:finite-path-upper}, while the lower
bound is obtained from a partition construction that applies to
every connected target. Thus, the function $b(x)$ introduced above
provides the natural coefficient in the resulting bounds.

\begin{theorem}\label{thm:path-dense-asymptotic}
Fix an integer $t\ge1$ and a real number $\beta>1$.
Suppose that $r,n\to\infty$, \( \frac rs\longrightarrow\beta, \ \log r=o(n)\). 
Then
\begin{equation}\label{eq:path-dense-sandwich}
 \bigl(b(\beta)^t-o(1)\bigr)n
 \le R(P_n;r,s,t)
 \le \bigl(b(\beta^t)+o(1)\bigr)n.
\end{equation}
The lower bound holds more generally with $P_n$ replaced by
any connected $n$-vertex graph.

In particular, if $\beta=m$ is an integer with $m\ge2$, then
\begin{equation}\label{eq:path-dense-exact}
 R(T_n;r,s,t)=(m^t+o(1))n
\end{equation}
for every $n$-vertex tree $T_n$.

When $t=1$, the condition $\log r=o(n)$ can be omitted, and
\begin{equation}\label{eq:path-dense-one-color}
 R(P_n;r,s,1)=\bigl(b(\beta)+o(1)\bigr)n.
\end{equation}
\end{theorem}

\begin{proof}
\emph{The upper bound.}
Recall $\alpha=\frac{\binom rt}{\binom st}$. 
Since $t$ is fixed and $r/s\to\beta$, we have $\alpha\to \beta^t$.
For $z\ge1$, define
\[
 F(z)=\lfloor z\rfloor+(z-\lfloor z\rfloor)^2.
\]
The function $z^2/F(z)$ is continuous and strictly increasing:
on an interval $j<z<j+1$, its derivative is
\[
 \frac{2zj(j+1-z)}{F(z)^2}>0.
\]
By the definition of $b$, we have $b(\beta^t)^2=\beta^t F(b(\beta^t))$.
Hence every fixed $z>b(\beta^t)$ satisfies \( z^2>\beta^t F(z)\). 
Take $N=\lceil zn\rceil$. The exact path extremal
formula~\eqref{eq:exact-path-ex} gives
\[
 \frac{2\operatorname{ex}(N,P_n)}{n^2}\longrightarrow F(z).
\]
Consequently, for all sufficiently large parameters,
\( 
 \binom N2>\alpha\,\operatorname{ex}(N,P_n)\).
Applying~\eqref{eq:finite-path-upper} and then letting
$z\downarrow b(\beta^t)$ proves the upper bound in
\eqref{eq:path-dense-sandwich}. 
We emphasize that the above upper-bound argument does not require the condition $\log r=o(n)$. 

\medskip
\noindent
\emph{The lower bound for connected graphs.}
Put $c=b(\beta)$, $\theta=1/c$, and $j=\lfloor\beta\rfloor$.
The definition of $b$ gives
\begin{equation}\label{eq:path-1-jtheta}
 j\theta^2+(1-j\theta)^2=\frac1\beta,
 \qquad
 \frac1{j+1}<\theta\le\frac1j.
\end{equation}
For every sufficiently small fixed $\varepsilon>0$, choose a
probability vector $(\pi_1,\ldots,\pi_L)$ such that
\[
 \max_a\pi_a\le\theta+\varepsilon<1,
 \qquad
 \rho:=\sum_{a=1}^L\pi_a^2>\frac1\beta.
\]
To see that this is possible when $\beta$ is not an integer,
start with $j$ entries equal to $\theta$ and one entry equal
to $1-j\theta$, and increase the first $j$ entries slightly
while decreasing the last one. The squared sum increases,
since its derivative at the original vector is
$2j((j+1)\theta-1)>0$.
When $\beta=j$ is an integer, start with the uniform
distribution on $j$ classes and replace two entries by
$1/j+\varepsilon$ and $1/j-\varepsilon$.

Choose \( h=\left\lfloor\min\{r,\sqrt{n\log r}\}\right\rfloor\). 
The assumptions imply \( h\to\infty,\ h=o(n),\ \log h=o(r),\ \log r=o(h)\). 
Write $\eta=(\theta+\varepsilon)^t$, fix
$0<\delta<\eta^{-1}-1$, and put \( N_0=\left\lfloor(\eta^{-1}-\delta)h\right\rfloor\). 
Independently for every $v\in V(K_{N_0})$ and $i\in[r]$,
choose a class $X_{v,i}\in[L]$ according to the probability
vector above. For each edge $uv$, let
\[
 I(uv)=\{i\in[r]:X_{u,i}=X_{v,i}\}.
\]
Then $|I(uv)|$ has distribution $\operatorname{Bin}(r,\rho)$.

For $T\in\binom{[r]}t$ and $\mathbf a\in[L]^T$, define
\[
 V(T,\mathbf a)
 =\{v:X_{v,i}=a_i\text{ for every }i\in T\}.
\]
Each vertex belongs to this set with probability at most
$\eta$. Hoeffding's inequality~\cite{Hoeffding1963} and the union bound show that
the probability that some edge has $|I(uv)|<s$ or some cell
has $|V(T,\mathbf a)|\ge h$ is at most
\[
 \binom{N_0}{2}
 \exp\!\left(-2r\left(\rho-\frac sr\right)^2\right)
 +
 \binom rt L^t
 \exp\!\left(-\frac{2(h-\eta N_0)^2}{N_0}\right).
\]
Here $\rho-s/r$ is bounded below by a positive constant,
and $h-\eta N_0\ge\eta\delta h$. The displayed expression
therefore tends to zero, because $\log h=o(r)$ and
$\log r=o(h)$.

Thus there is a realization for which every edge has at
least $s$ available colors and every cell has fewer than
$h$ vertices. Choose an arbitrary $s$-subset
$\phi(uv)\subseteq I(uv)$ for each edge.
For every $T\in\binom{[r]}t$, each component of the
$T$-support lies in a single cell $V(T,\mathbf a)$ and
therefore has at most $h-1$ vertices.

Now replace each vertex by a set of \( q=\left\lfloor\frac{n-1}{h-1}\right\rfloor \) 
vertices. Give every edge between two such sets the label
of the corresponding original edge, and assign arbitrary
$s$-element labels to edges within each set.
Every component of a new $T$-support is contained in the
union of the vertex sets corresponding to one original
$T$-support component. Its order is therefore at most 
\( q(h-1)\le n-1\). 
The resulting coloring contains no $t$-monochromatic copy
of any connected $n$-vertex graph.

Since $h\to\infty$ and $h=o(n)$, we have $qh/n\to1$.
Thus, for every connected $n$-vertex graph $G_n$,
\[
 R(G_n;r,s,t)>qN_0
 =\bigl((\theta+\varepsilon)^{-t}-\delta+o(1)\bigr)n.
\]
Take the joint limit with $\varepsilon,\delta$ fixed,
and then let $\delta\downarrow0$ and
$\varepsilon\downarrow0$. This gives
\[
 R(G_n;r,s,t)
 \ge\bigl(\theta^{-t}-o(1)\bigr)n
 =\bigl(b(\beta)^t-o(1)\bigr)n,
\]
proving the general lower bound and, in particular,
\eqref{eq:path-dense-sandwich}.

\medskip

\medskip
\noindent
Suppose that $\beta=m\in\mathbb Z_{\ge2}$.
For every $n$-vertex tree $T_n$, since $b(\beta)^t=m^t$ and $\alpha\to m^t$, 
the lower bound just proved and the tree bound~\eqref{eq:tree-bound} give
\[
 (m^t-o(1))n
 \le R(T_n;r,s,t)
 \le \lfloor\alpha(n-2)\rfloor+2
 =(m^t+o(1))n.
\]
This proves~\eqref{eq:path-dense-exact}.

\medskip
\medskip
\noindent
For the case $t=1$,  note that the condition $\log r=o(n)$ is not needed. Fix $1<c_0<b(\beta)$.
Choose the probability vector above with
$\theta+\varepsilon<1/c_0$.
By approximating its entries by proportions of a
sufficiently large fixed integer $M$, choose a partition
of $[M]$ with block sizes $k_1,\ldots,k_L$ satisfying
\[
 p_*:=\frac{\sum_{a=1}^L\binom{k_a}{2}}{\binom M2}
 >\frac1\beta,
 \qquad
 \frac{\max_a k_a}{M}<\frac1{c_0}.
\]

Take the list of all $M!$ permutations of this partition,
with multiplicity. Every pair of distinct elements of $[M]$
lies in the same block in exactly $p_*M!$ members of the list.
Repeat the list to obtain $r$ partitions, filling any
remaining positions with partitions from the same list.
Every pair of distinct elements then lies in a common
block in at least
\[
 p_*M!\left\lfloor\frac r{M!}\right\rfloor
 \ge p_*r-M!
 \ge s
\]
of these partitions, for all sufficiently large $r$.

Replace the $M$ elements by nearly equal vertex sets of
total order $N=\lfloor c_0n\rfloor$.
For each of the $r$ partitions, form a color support by
making the union of the vertex sets in each block a clique.
Every edge between different vertex sets belongs to at
least $s$ supports, while every edge within a vertex set
belongs to all $r$ supports. Retain any $s$ available
colors on each edge.

Every component of each color support has order at most
\( \bigl(\max_a k_a\bigr)
 \left\lceil\frac NM\right\rceil<n \) 
for all sufficiently large $n$.
Hence this coloring avoids every connected $n$-vertex
graph, and in particular
\[
 R(P_n;r,s,1)>\lfloor c_0n\rfloor.
\]
The partition and $M$ are fixed before $r,n\to\infty$,
so no restriction on their relative growth is needed.
Letting $c_0\uparrow b(\beta)$ and combining this lower
bound with the upper bound proves
\eqref{eq:path-dense-one-color}.
\end{proof}

\section{Complete graphs}
\label{sec:cliques}

Recall that $\alpha=\binom rt/\binom st$, and write
$R(n;r,s,t)=R(K_n;r,s,t)$. We first give general upper
and lower bounds for $R(n;r,s,t)$.

\begin{theorem}\label{thm:KnlowUp}
There exist absolute constants $c_1,C_1>0$ such that,
for every fixed choice of integers $r>s\ge t\ge1$,
there is a positive constant $\eta=\eta(r,t)$ for which
the following bounds hold for all sufficiently large $n$:
\begin{equation}\label{eq:clique-dense}
 \eta\,2^{c_1n\max\left\{
       \frac{r-s}{r-t+1},
       \frac{(r-s)^3}{(r-t+1)^2}
       \right\}}\le 
 R(n;r,s,t)
 \le
 2^{C_1n\binom rt
       \left(1-\frac1\alpha\right)^2
       \log_2\frac{1}{
         \min\left\{\frac1\alpha,\,1-\frac1\alpha\right\}}}.
\end{equation}
Moreover, if $\lfloor r/t\rfloor>r-s$, then the exponent
in the upper bound can be replaced by
\[
 C_1n\min\left\{
 \binom rt
   (1-\frac1\alpha)^2
   \log_2\frac{1}{
     \min\left\{\frac1\alpha,\,1-\frac1\alpha\right\}}, 
 \frac{(r-s)^2}{\lfloor r/t\rfloor}
   \log_2
   \frac{\lfloor r/t\rfloor}{
     \min\left\{\lfloor r/t\rfloor-r+s,\,r-s\right\}}
 \right\}.
\]
\end{theorem}

\begin{proof}
We first recall the lower bound of Conlon et
al.~\cite[Theorem~1.1]{CFHMSV}: for fixed integers
$r>s\ge1$ and all sufficiently large $n$,
\[
 R(n;r,s,1)
 \ge
 \eta(r,1)\,
 2^{c_1n\max\left\{
       \frac{r-s}{r},
       \frac{(r-s)^3}{r^2}
       \right\}},
\]
where $c_1>0$ is absolute and $\eta(r,1)>0$ depends
only on $r$.
By Proposition~\ref{prop:padding},
\[
 R(n;r,s,t)\ge R(n;r-t+1,s-t+1,1).
\]
Applying the preceding lower bound with $r,s$ replaced
by $r-t+1,s-t+1$, respectively, gives
\[
 R(n;r,s,t)
 \ge
 \eta\,
 2^{c_1n\max\left\{
       \frac{r-s}{r-t+1},
       \frac{(r-s)^3}{(r-t+1)^2}
       \right\}},
\]
where $\eta=\eta(r,t)>0$ is the prefactor corresponding
to the palette size $r-t+1$ and therefore depends only
on $r,t$. This proves the lower bound
in~\eqref{eq:clique-dense}.

For the upper bound, the same theorem of Conlon et al.
gives, for fixed $r>s\ge1$ and all sufficiently large $n$,
\[
 R(n;r,s,1)
 \le
 2^{C_1n\frac{(r-s)^2}{r}
       \log_2\frac{r}{\min\{s,r-s\}}},
\]
where $C_1>0$ is absolute.
The shadow comparison~\eqref{shadow1} gives
\[
 R(n;r,s,t)
 \le R\left(n;\binom rt,\binom st,1\right).
\]
Applying the preceding upper bound with $r,s$ replaced
by $\binom rt,\binom st$, respectively, and using
$\binom st/\binom rt=1/\alpha$, we obtain
\[
 R(n;r,s,t)
 \le
 2^{C_1n\binom rt
       \left(1-\frac1\alpha\right)^2
       \log_2\frac{1}{
         \min\left\{\frac1\alpha,\,1-\frac1\alpha\right\}}}.
\]
This proves the upper bound in~\eqref{eq:clique-dense}.

Finally, suppose that $\lfloor r/t\rfloor>r-s$.
By disjoint-witness reduction~\eqref{eq:block},
\[
 R(n;r,s,t)
 \le
 R\left(
 n;\left\lfloor\frac rt\right\rfloor,
 \left\lfloor\frac rt\right\rfloor-r+s,1
 \right).
\]
The two color parameters on the right differ by $r-s$.
Applying the upper bound of Conlon et al. once more gives
\[
 R(n;r,s,t)
 \le
 2^{C_1n
   \frac{(r-s)^2}{\lfloor r/t\rfloor}
   \log_2
   \frac{\lfloor r/t\rfloor}{
     \min\left\{\lfloor r/t\rfloor-r+s,\,r-s\right\}}}.
\]
Taking the minimum of this bound and the shadow upper
bound proves the final assertion. Since all color
parameters in these applications are fixed, the bounds
hold simultaneously for all sufficiently large $n$.
\end{proof}

Now fix the clique order $n$ and let the palette grow. The
relevant quantities are the proportion $s/r$ of colors assigned
to each edge and the proportion $t/r$ required on the whole clique.

Define $\Phi_n:[0,1]\to[0,1]$ by linearly interpolating the points
\[
 \left(1-\frac1j,\frac{(j)_n}{j^n}\right),\qquad j=1,2,\ldots,
 \qquad \Phi_n(1)=1,
\]
where $(j)_n=j(j-1)\cdots(j-n+1)$, with $(j)_n=0$ for $j<n$.
These coordinates are the limiting edge and clique densities of balanced
complete $j$-partite graphs. 
Its successive slopes are nondecreasing, so $\Phi_n$ is convex. 
In particular, $\Phi_n(x)=0$ for $x\le1-1/(n-1)$.  

\begin{theorem}\label{thm:clique-transition}
Fix $n\ge3$ and $0<\tau\le\sigma<1$. Suppose that, as
$r\to\infty$, the integers $r\ge s\ge t\ge1$ satisfy
\(
 \frac sr\longrightarrow\sigma\),
 and  \(\frac tr\longrightarrow\tau\).
Then
\[
 \begin{cases}
 R(K_n;r,s,t)=O(1),&\text{if }\tau<\Phi_n(\sigma),\\[2pt]
 R(K_n;r,s,t)\ge\exp(cr),&\text{if }\tau>\Phi_n(\sigma),
 \end{cases}
\]
where $c>0$ and the constant implicit in $O(1)$ depend only on
$n,\sigma,\tau$.
\end{theorem}

For example, $\Phi_3(x)=\frac43x-\frac23$ on $[1/2,2/3]$.
Thus, when $s/r\to3/5$, the threshold for triangles is $2/15$:
the Ramsey numbers stay bounded if $t/r\to\tau<2/15$, and grow
at least exponentially in $r$ if $t/r\to\tau>2/15$.

\begin{proof}
\emph{The bounded regime.}
Suppose $\tau<\Phi_n(\sigma)$. By continuity, choose a fixed
$N\ge n$ such that
\[
 \Phi_n\bigl(\sigma(1-1/N)\bigr)>\tau.
\]
Consider an $(r,s)$-set-coloring of $K_N$ with no
$t$-monochromatic $K_n$. For each $i\in[r]$, let $J_i$ consist
of the edges whose labels contain $i$. Every edge belongs to
exactly $s$ of these graphs, while every $n$-set spans a clique
in at most $t-1$ of them. Hence
\[
 \frac1r\sum_{i=1}^r\frac{2e(J_i)}{N^2}
 =\frac sr(1-1/N),\qquad
 \sum_{i=1}^r\#K_n(J_i)\le(t-1)\binom Nn,
\]
where $\#K_n(H)$ denotes the number of unlabelled $n$-cliques in $H$.
Bollob\'as's piecewise-linear lower bound for clique densities~\cite{Bollobas,LMR} gives
\begin{equation}\label{eq:bollobas}
 \frac{n!\,\#K_n(H)}{N^n}
 \ge\Phi_n\left(\frac{2e(H)}{N^2}\right)
 \qquad(|V(H)|=N).
\end{equation}
This normalization follows from the asymptotic form by taking
independent blow-ups of $H$. Since $\Phi_n$ is convex, averaging
\eqref{eq:bollobas} over the color supports yields
\begin{equation}\label{eq:clique-jensen}
 \Phi_n\left(\frac sr(1-1/N)\right)
 \le\frac{n!}{rN^n}\sum_{i=1}^r\#K_n(J_i)
 \le\frac{t-1}{r}.
\end{equation}
As $r\to\infty$, this contradicts the choice of $N$. Therefore
$R(K_n;r,s,t)\le N$ for all sufficiently large $r$.

\emph{The exponential regime.}
Suppose $\tau>\Phi_n(\sigma)$. Choose $p\in(\sigma,1)$ such that
$\Phi_n(p)<\tau$. There are an integer $j\ge1$ and
$\lambda\in[0,1]$ with
\[
 p=(1-\lambda)\left(1-\frac1j\right)
   +\lambda\left(1-\frac1{j+1}\right).
\]
On a vertex set of order $N$, construct a random graph as follows.
Choose $a=j$ with probability $1-\lambda$ and $a=j+1$ with
probability $\lambda$. Assign each vertex independently and
uniformly to one of $a$ classes, and join vertices in different
classes. A fixed pair is an edge with probability $p$, and a fixed
$n$-set is a clique with probability
\[
 (1-\lambda)\frac{(j)_n}{j^n}
 +\lambda\frac{(j+1)_n}{(j+1)^n}=\Phi_n(p).
\]
Take $r$ independent graphs $J_1,\ldots,J_r$ from this model.
For a pair $e$, let $X_e$ count the graphs containing it; for a
$n$-set $U$, let $Y_U$ count the graphs in which it is a clique.
Then
\[
 X_e\sim\operatorname{Bin}(r,p),\qquad
 Y_U\sim\operatorname{Bin}(r,\Phi_n(p)).
\]
Choose
\[
 0<c<\min\left\{(p-\sigma)^2,
                  \frac{2}{n}\bigl(\tau-\Phi_n(p)\bigr)^2\right\}
 \quad\text{and set}\quad N=\lfloor e^{cr}\rfloor.
\]
Hoeffding's inequality and a union bound show that the probability
of some $X_e<s$ or some $Y_U\ge t$ is at most
\[
 \binom N2 e^{-2r(p-s/r)^2}
 +\binom Nn e^{-2r(t/r-\Phi_n(p))^2},
\]
which tends to zero. Choose an outcome with $X_e\ge s$ for every
edge and $Y_U<t$ for every $n$-set. Label each edge $e$ by any
$s$ members of $\{i:e\in E(J_i)\}$. These labels have fewer than
$t$ common colors on every $K_n$, so $R(K_n;r,s,t)>N$.
\end{proof}

The construction realizes each segment of $\Phi_n$ by mixing two
adjacent multipartite models. For triangles at density $3/5$,
independent random edges would instead give the larger clique
probability $(3/5)^3=27/125$; the mixture gives the threshold $2/15$.

\section{The intersection spectrum and further questions}
\label{sec:questions}

For fixed $G,r,s$, the sequence
\begin{equation}\label{eq:spectrum}
 R(G;r,s,1)\le R(G;r,s,2)\le\cdots\le R(G;r,s,s)
 =R\left(G;\binom rs\right)
\end{equation}
connects set-coloring Ramsey theory with ordinary multicolor
Ramsey theory. The spectrum reveals how requiring more common colors
changes both extremal structure and growth rate.

\begin{question}
For which graphs and palettes are the inequalities in
\eqref{eq:spectrum} strict? How large can consecutive jumps be?
\end{question}

Even the dense clique slice is instructive.
When \(s=r-1\) and \(t\le r/2\),
Theorem~\ref{thm:KnlowUp} gives two-sided exponential
bounds, with an improved upper exponent from the
disjoint-witness reduction. At \(t=r-1\), the parameter
reduces to the ordinary \(r\)-color Ramsey number.
Determining the transition between these scales would
extend the questions underlying~\cite{ES,CFHMSV,ACMMM}.
Determining the growth on the critical curve
$\tau=\Phi_n(\sigma)$ in Theorem~\ref{thm:clique-transition}
is another natural problem.

\begin{question}
When does equality hold in the inequality
\eqref{shadow1}, and how large can its 
additive or multiplicative gap be?
\end{question}
For stars, Corollary~\ref{cor:star-design} gives equality for
infinitely many orders in design regimes, whereas
Section~\ref{sec:separation} gives a strict gap already for $P_3$.
For fixed $s,t$ and growing $r$, our packing bounds give
asymptotic equality with the shadow problem for long paths and,
in the stated iterated limit, even cycles. For a fixed palette,
obtaining the optimal coefficients remains open.

Theorem~\ref{thm:star-packing} leaves a concrete realization problem. 
When \(\nu(r,s,t;n-1)\) is even but no maximum multiset has every label occurring with even multiplicity,
can a coloring of the complete graph \(K_{\nu(r,s,t;n-1)+1}\) realize maximum local packings at every vertex? 
Under additional divisibility conditions, 
Theorem~\ref{thm:rigid-star} provides a partial answer through an incidence-rigidity argument, 
showing that such a coloring forces every label class to be a regular spanning subgraph. 
This rigidity yields exact star Ramsey numbers in the specified parameter regime. 
A characterization beyond this regime may lead to further exact star cases and shed light on 
how local packing structures constrain global set-colorings.

\section*{Disclosure of AI use}

This paper were completed approximately one year ago, with the initial draft essentially finalized by late last year. 
This year, we engaged in a revision process to refine and supplement the manuscript. 
During this revision phase, we utilized AI tools in specific aspects. 
First, we asked the AI to search for explicit example, 
which led to the parameter values presented in Section 2. 
Second, we used AI to polish the English and verify the logic of certain proofs. 
Beyond these, after we had improved the upper bounds for paths, established the packing lower bounds, 
we directed the AI to investigate dense-palette regime and it suggested an approximate lower bound for the case where s is close to r. 
By combining our upper bounds with this AI-suggested asymptotic lower bound, 
we established Dense-palette asymptotics (Theorem 5.6).
We have thoroughly reviewed and verified AI-assisted modifications. 
All core mathematical ideas and conclusions are our own. 
We take full responsibility for all mathematical claims, proofs, and the final content of the paper.

\bibliographystyle{plainnat}
\bibliography{references}
\end{document}